\documentclass{article}

\usepackage[naustrian, english]{babel} 
\usepackage[T1]{fontenc}
\usepackage{lmodern}
\usepackage{textcomp}
\usepackage{amssymb} 
\usepackage{amsthm}
\usepackage{mathtools}
\usepackage[pdftex]{graphicx}
\usepackage[export]{adjustbox}
\usepackage{caption}
\usepackage{color}
\usepackage[arrow, matrix, curve]{xy}
\usepackage{lmodern}
\usepackage[utf8]{inputenc}
\usepackage{CJKutf8}
\usepackage{url}
\usepackage{authblk}
\usepackage{hyperref}
\usepackage{float}
\usepackage{listings}
\usepackage{appendix}
\usepackage{pdfpages}
\usepackage{algorithm}
\usepackage{algorithmic}
\usepackage{geometry}
\usepackage{bbm}
\usepackage{bbold}
\usepackage{graphicx}
\usepackage{wrapfig}
\usepackage{caption}
\usepackage{subcaption}
\title{A network-level transport model of tau progression in the Alzheimer's brain}
\author{Veronica Tora \thanks{Department of Mathematics, University of Roma ``Tor Vergata'', Roma, Italy {\tt  veronica.tora2@unibo.it, Veronica.Tora@uniroma2.it}} \and  
Justin Torok \thanks{Department of Computational Biology and Medicine, Weill Cornell Medical School, New York, New York, United States of America {\tt Justin.Torok@ucsf.edu}} \and
Michiel Bertsch \thanks{Department of Mathematics, University of Roma ``Tor Vergata'', Roma, Italy, Istituto per le Applicazioni del Calcolo ``M. Picone'', Consiglio Nazionale delle Ricerche, Roma, Italy {\tt bertsch@mat.uniroma2.it}}\and 
Ashish Raj \thanks{Department of Radiology and Biomedical Imaging, University of California at San Francisco, San Francisco, California, United States of America {\tt ashish.raj@ucsf.edu}} \and }

\date{}

\begin{document}

\maketitle

\begin{abstract}
One of the hallmarks of Alzheimer's disease (AD) is the accumulation and spread of toxic aggregates of tau protein. 
The progression of AD tau pathology is thought to be highly stereotyped, which is in part due to the fact that tau can spread between regions via the white matter tracts that connect them. 
Mathematically, this phenomenon has been described using models of ``network diffusion,'' where the rate of spread of tau between brain regions is proportional to its concentration gradient and the amount of white matter between them. 
Although these models can robustly predict the progression of pathology in a wide variety of neurodegenerative diseases, including AD, an underexplored aspect of tau spreading is that it is governed not simply by diffusion but also active transport along axonal microtubules.
Spread can therefore take on a directional bias, resulting in distinct patterns of deposition, but current models struggle to capture this phenomenon.
Recently, we have developed a mathematical model of the axonal transport of toxic tau proteins that takes into account the effects tau exerts on the molecular motors. 
Here we describe and implement a macroscopic version of this model, which we call the Network Transport Model (NTM). 
A key feature of this model is that, while it predicts tau dynamics at a regional level, it is parameterized in terms of only microscopic processes such as aggregation and transport rates; that is, differences in brain-wide tau progression can be explained by its microscopic properties. 
We  provide numerical evidence %demonstrate 
that, as with the two-neuron model that the NTM extends, there are distinct and rich dynamics with respect to the overall rate of spread and the staging of pathology when we simulated the NTM on the hippocampal subnetwork. 
The theoretical insights provided by the NTM have broad implications for understanding AD pathophysiology more generally.

\end{abstract}
\medskip

\noindent{\bf Keywords:} Alzheimer's disease, models on graphs, Tau protein, transport-reactions equations,
numerical simulations.

\clearpage
%
%
%
%
%

%\noindent \jtcomment{This is what I'm envisioning more specifically for the Introduction:}

\section{Introduction} 
%\begin{enumerate}
%    \item \textit{2-3 paragraphs of biological motivation}
%    \item \textit{Previous network-level models}
%    \item \textit{Directionality as an underexplored mechanism}
%    \item \textit{Brief introduction to the single-axon model from Torok, et al.}
%\end{enumerate}

% ******************** AR Portion Begins ***************

Alzheimer's disease (AD) is a neurodegenerative disorder whose key pathological hallmark is the abnormal deposition of microtubule-associated protein tau and its progressive ramification throughout the brain in a characteristic and highly stereotyped spatiotemporal pattern. 
Tau tangles appear first in the locus coeruleus and then spread to the entorhinal cortex, hippocampus, temporal areas, and finally throughout the cortex~\cite{Braak1991,Braak2006d}. 
Significant \textit{in vitro} and \textit{in vivo} evidence exists that tau predominantly migrates \textit{trans-synaptically}, while white matter tracts between regions serve as the conduits for the transmission of tau from affected regions to unaffected regions ~\cite{Clavaguera2009,Katsikoudi2020,Liu2012,Wang2017a}. 

Understanding the mechanism behind how the disease process unfolds over time is critical for ultimately finding effective treatments for AD. 
Mathematical modeling of tau propagation can provide a platform for integrating experimental data into a cohesive theoretical framework that may be used to test hypotheses for which direct evidence is difficult to acquire. 
The key insight for modeling trans-synaptic spread of tau is that it can be approximated by a \textit{graph diffusion model}, where discrete gray-matter regions are the vertices of a graph and the structural connectivity values for all region pairs are its edge weights. 
More specifically, the Network Diffusion Model (NDM)~\cite{Raj2012} and subsequent connectome-based spread models~\cite{Fornari2019,Iturria-Medina2014,Schafer2020,Weickenmeier2018,Weickenmeier2019,Bertsch2023} posit that, given an initial distribution of pathology in the brain, the regional pathology at future time points is a function of the concentration differences and connectivity values between all region pairs. 
Remarkably, despite the complexity of AD at a molecular and cellular level, these simple, macroscopic models recapitulate the canonical Braak staging of AD~\cite{Weickenmeier2018,Weickenmeier2019} as well as pathology progression in human subjects~\cite{Raj2015,Schafer2020,Raj2021a}. 
%Not only that, but other disorders featuring the trans-synaptic spread of pathological protein species, such as Parkinson's disease (PD)~\cite{Pandya2019,Yau2018} and amyotrophic lateral sclerosis (ALS)~\cite{Pandya2021}, have also been modeled as a graph diffusion process with similarly strong results. 

Despite the success of connectome-based spread models, they do not have sufficient complexity to capture non-passive diffusive spread. 
Two key aspects remain unaddressed in current mathematical models: 
\begin{enumerate}
    \item \textbf{Oligomer conversion and aggregation kinetics}. 
    Connectome models of tau spread have historically limited themselves to modeling the net effect of these proteins, regardless of their oligomeric diversity. 
    In fact, it is well known that tau oligomers of varying sizes engage in kinetic equilibrium, and these aggregation and fragmentation processes strongly contribute to their ability to spread throughout the brain~\cite{Kundel2018,Xue2009}. 
    The modeling of the kinetics of protein aggregation are rather well established by us and others ~\cite{Bertsch2016, Franchi2012, Franchi2016, Franchi2017, Smoluchowski1918}. 
    However, its incorporation within a network spread model has received far less attention. 

    \item \textbf{From passive diffusion to active transport}. 
    Tau may not migrate in a purely \textit{diffusive} manner along concentration gradients, as is commonly assumed in connectome-based spread models, but instead undergoes active transport via molecular motors attached to microtubules in either \textit{anterograde} (i.e., with axon polarity) or \textit{retrograde} (i.e., against axon polarity) directions. 
    The hyperphosphorylation of pathological tau disrupts its direct interactions with microtubules~\cite{Alonso1994,Alonso1997,Konzack2007} and the motor proteins themselves~\cite{Cuchillo-Ibanez2008,Rodriguez-Martin2013,Stern2017}. 
    Together, these effects lead to aberrant axonal transport and the mis-sorting of tau into the neuronal somatodendritic compartment~\cite{Ballatore2007,Zempel2017}. 
    To some extent active transport along fiber projections may be approximated by anisotropic diffusion, an approach taken in recent studies \cite{Kuhl2019,Henderson2019}. 
    However, no diffusive process can fully capture transport processes, and in particular a key aspect of transport: directionality of transmission. 
    It was noted by a recent review article that the emergence of net retrograde bias in tau propagation may be explainable via its interaction with transport kinetics \cite{Torok2021}.   
    
\end{enumerate}

It is therefore critical to model active transport of tau on fiber projections, along with network connectivity. 
At a microscopic level, Torok, \textit{et al.} recently proposed a two-species, multicompartment model to explore the interactions between the pathological axonal transport of tau and the formation and breakdown of insoluble tau aggregates~\cite{Torok2021}. 
Leveraging insights from \textit{in vitro} work demonstrating that the primary anterograde-directed motor protein, kinesin-1, has increased activity in the presence of monomeric hyperphosphorylated tau \cite{Cuchillo-Ibanez2008,Rodriguez-Martin2013,Stern2017} and is knocked down by tau aggregates~\cite{Sherman2016}, the mathematical model poses a simple tau-concentration-dependent feedback mechanism on the anterograde velocity of tau transport. 
The authors interrogated the complex dynamics that emerge between the interplay of tau aggregation and axonal transport feedback, finding that higher aggregation rates generally led to stronger retrograde biases in tau deposition at steady-state. 
This work represents one of the first attempts to connect the microscopic properties of tau conformers to pathological changes that can be observed at macroscopic timescales. 

\begin{figure}
\centering
\includegraphics[width=0.75\linewidth]{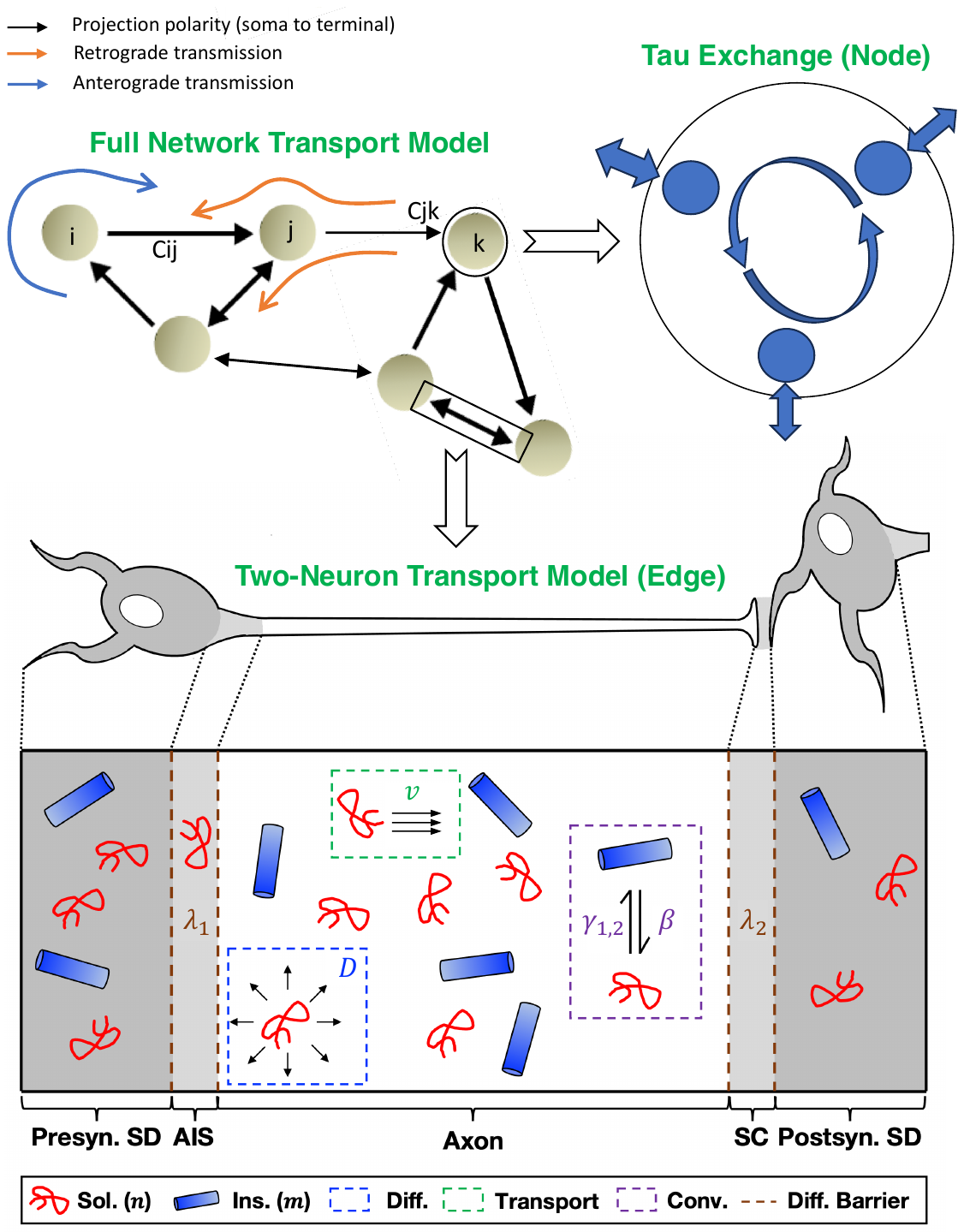}
\caption{\textbf{Illustration of the model system.} 
At the whole network level, brain regions are represented by nodes and white matter fiber projections between them by edges (top panel). 
Tau pathology propagates on this network in an anterograde or retrograde direction, depending on the cell polarity and the properties of tau itself. 
Instead of passive graph diffusion previously used to model the transmission along an edge, here we use the active axonal transport model from Torok, \textit{et al.}, which is schematized version in the bottom panel. 
It models two distinct species of pathological tau, soluble (red) and insoluble (blue), within a multi-compartment, two-neuron system mimicking the single-edge system shown in the top panel. 
The main biological phenomena captured in this model are diffusion (blue box), active transport (green box), species interconversion through fragmentation and aggregation (purple box), and a diffusion-based barrier to inter-compartmental spread (brown dashed lines).
The full model involves iteratively solving the PDE-based, single-edge model for the fluxes at the boundaries of the system and the resulting solving the mass exchange problem at the nodes (right panel).
Sol. - soluble tau; Ins. - insoluble tau; Diff. - diffusion; Conv. - tau interconversion; Diff. Barrier - diffusion barrier; Presyn. SD - presynaptic somatodendritic compartment; AIS - axon initial segment; SC - synaptic cleft; Postsyn. SD - postsynaptic somatodendritic compartment.
Figure adapted from the original manuscript~\cite{Torok2021}. 
% edge model therefore determines the mass fluxes, which are considered to engage in an exchange process within the node (right panel). 
% The full model evolution therefore involves successively solving the PDE at each edge and updating the exchange system within each node. 
}
\label{fig:schematic}
\end{figure}

Here we provide a macroscopic model which combines the dynamics of soluble and insoluble tau in the gray-matter regions and 
Torok, \textit{et al.} axonal transport model~\cite{Torok2021} in the white-matter tracts. 
This model, which we call the Network Transport Model (NTM), enables us to simulate the dynamics of soluble and insoluble tau in terms of the diffusion-advection and aggregation-fragmentation processes as in~\cite{Torok2021}, but at the network level. 
The essential elements of the NTM are illustrated in \textbf{Figure~\ref{fig:schematic}}.

The full network model involves a transport-reaction PDE on each edge, where the dynamics of soluble and insoluble tau within the white matter tracts is governed by the Torok, \textit{et al.} model, and a diffusion-reaction equation on the nodes describing tau dynamics in the gray matter regions. 
A straightforward mass transfer mechanism of soluble tau between edges and nodes determines the incoming flux of soluble tau into the nodes. 
The NTM therefore provides a necessary augmentation of prior network-based models of tauopathy, such as those explored by our group~\cite{Raj2012,Raj2015,Raj2021a,Bertsch2023} and others~\cite{Iturria-Medina2013,Weickenmeier2018,Henderson2019}, as it allows for the characterization of macroscopic tauopathy dynamics in terms of the microscopic properties of soluble and insoluble tau species, and their active transport along axons. 

While the full NTM is computationally infeasible to simulate on the full network, we provide and implement a quasi-static approximation to the NTM that maintains the basic properties of the full NTM and is more tractable numerically.
We propose in this paper some key mathematical innovations necessary to make this approach feasible.
First, we show that the full, connectome-coupled system of transport PDEs can be reasonably approximated by separating the dynamics into two time scales: a fast time scale (order of hours to days) whereby transport processes in individual edges are established; followed by a slow time scale (order of months to years) over which the full network couples to these individual egde processes via slow exchange within network nodes. 
These two regimes are schematically illustrated in \textbf{Figure~\ref{fig:schematic}}. 
Second, by assuming this separation of time scales, the dynamics of the ``fast'' processes at the terminals of each edge can be approximated by a steady-state solution to the PDE-based model, which we compute using numerical integration. 
The major modeling challenge is to determine the correct mathematical description of the mass transfer between edges and nodes at the ``slow'' time scale. 
The resulting local mass balance between a node and its incoming and outgoing edges contains, in addition to mass fluxes entering the node, a feedback mechanism. 

We demonstrate that the resulting quasi-static NTM model exhibits rich dynamics as a function of these microscopic parameters, expanding the range of behaviors that can be exhibited with previous models, particularly in terms of directionally biased flows on the connectome. 
The quasi-static NTM is therefore capable of delivering more robust insights into how specific tau species may differentially propagate in the brain, which has important clinical implications for our understanding and treatment of AD and other tauopathic dementias.
%

% ******************** AR Portion ends ***************

%
%
%
\section{Modeling}
%
%
%
%
%
% \mbcomment{I prefer to skip the first part of this section; in Section 2.4.1 there is now a first attempt to reassume the quasi-static NTM model}

Here we describe the development of the Network Transport Model.  
In \textbf{Section 2.1}, we describe in detail what constitutes the connectivity graph, or connectome, on which we simulate tau spreading, both in terms of its mathematical formulation and the empirical data on which it is based. 
We then give a high-level overview of the original two-neuron axonal transport model from Torok, \textit{et al.} \cite{Torok2021}, as well as several of the modifications required for using this model to describe the internal dynamics of tau within each edge on the network (\textbf{Section 2.2}). 
For further details, we refer the reader to the original publication.
The full NTM in the continuous case is proposed in \textbf{Section 2.3}, which incorporates the axonal transport model of tau within the white matter tracts (edges of the connectivity graph) with the dynamics of tau in the gray matter regions (nodes of the connectivity graph). 
The latter is described by a diffusion-reaction equation on the nodes, where the diffusion mechanism is determined by the contributions of the mass fluxes along the edges.
Lastly, by assuming a separation of time scales for the tau dynamics on the edges and tau dynamics in the nodes, we derive the quasi-static approximation to the NTM, where a local mass balance problem has to be solved to determine the mass exchange of tau between a node and its incoming and outgoing edges, (\textbf{Section 2.4}).
\subsection{The structural connectivity graph}
The ``wiring diagram'' of the brain can be described in terms of a structural connectivity graph, $G$, where edges represent white matter tracts and vertices represent gray matter regions.  
More precisely, we define $G$ to be a weighted, directed graph with a finite number $h$ of vertices $P_i$ and edges $e_{ij}$ ($i\ne j$) directed from vertex $P_i$ to $P_j$. 
We distinguish edges $e_{ij}$ and $e_{ji}$ by their polarization: the polarization of $e_{ij}$ is directed from $P_i$ to $P_j$, whereas $e_{ji}$ is directed from $P_j$ to $P_i$. 

The graph is endowed with a weight function $c$ such that:
\begin{equation}\label{weights}
c(P_i,P_j)=
\begin{cases}
c_{ij}>0 & \text{if $P_i$ and $P_j$ are adjacent}\;. \\
0 & \text{otherwise}
\end{cases}
\end{equation}
The connectivity weights $c_{ij}$ express the strength of the connection between the $i$-th and $j$-th brain compartment. In a directed graph the weight function is not required to be symmetric; therefore, the weight $c_{ij}$ can be different from $c_{ji}$.
 
Here we utilize the mouse mesoscale connectivity atlas (MCA) from the Allen Institute for Brain Science~\cite{Oh2014}, which uses viral tracing methods to determine both the weights and polarity of brain connections at a fine regional parcellation. 
We simulate the NTM on the hippocampal subcircuit of the connectome, which has particular relevance to AD: this includes 11 hippocampal structures, 3 retrosplenial areas from the neocortex, and the piriform area, for a total of 30 regions across both hemispheres.
\subsection{A single-edge model}\label{S edges}
%
%
%
%
%
%\jtcomment{I changed rho to phi here and throughout.}
A mathematical description of the axonal transport dynamics of pathological tau was originally developed by Torok, \textit{et al.}, which postulates that the spreading of tau between neurons depends not only on the concentration difference between them, as has been previously proposed \cite{Raj2012,Iturria-Medina2013,Weickenmeier2018}, but also on the interactions between tau and the molecular motors of the axon \cite{Kuznetsov2016,Kuznetsov2018,Torok2021}. 
Its novelty lies in the fact that it can describe \textit{directionally biased spreading} of pathological tau; that is, preferential migration from presynaptic to postsynaptic neurons (anterograde) or from postsynaptic to presynaptic regions (retrograde), a process that has been described for tau \cite{Mezias2020a} and other prion-like proteins \cite{Henderson2019}. 
This model provides the basis of the dynamics on edges of the connectivity graph, as we detail below.
Here we make the key assumption that the edges of the connectivity graph (i.e., the white matter tracts of the brain), can be represented by bundles of connected neurons with \textit{independent tau dynamics}. 
Therefore, the fluxes at each vertex-edge boundary are equivalent to the single-axon fluxes multiplied by a proportionality constant given by the edge weights of the connectivity graph (see \textbf{Section 2.1}). 
For all definitions for the terms of this model, we refer the reader to \textbf{Table~\ref{Table_for_Variables}}.

Following Torok, \textit{et al.}, for a single bundle of neurons within edge $e_{ij}$, we describe the position within the two-neuron system by a 1D variable $x\in [0,L]$, where $L$ is the total ``length'' of the system. 
Because of the axon's high aspect ratio and the fact that the microtubules along which tau is transported are aligned with the long axis of the axon, we assume the dynamics occur predominantly in one dimension only.

We distinguish five segments, which represent biological compartments with distinct tau dynamics:
\begin{itemize}
\item[$(i)$] {\bf Presyn.\ SD}: presynaptic somatodendritic compartment, $(0,x_1)$. 
\item[$(ii)$] {\bf AIS}: axon initial segment, $(x_1, x_2)$. 
\item[$(iii)$] {\bf Axon}: axonal component, $(x_2, x_3)$. 
\item[$(iv)$] {\bf SC}: synaptic cleft, $(x_3, x_4)$. 
\item[$(v)$] {\bf Postsyn.\ SD}: postsynaptic somatodendritic compartment, $(x_4, L)$. 
\end{itemize}
We use the same compartment sizes as previously described and, for simplicity, assume that $L$ is a constant for all $e_{ij}$. 
 
Let $m_{ij}(x,t)$ and $n_{ij}(x,t)$ denote the densities at time $t$ per unit volume of insoluble and
soluble pathological tau respectively, where $t$ refers to the slow time scale. 
We introduce the following expressions for diffusive flux ($j_{\rm diff}$) and transport flux ($j_{\rm active}$) of $n=n_{ij}$:
$$
j_{\rm diff}(n_x)=-Dn_x, \qquad  j_{\rm active} (m,n)= v(m,n)n
$$
where $D>0$ is the diffusivity of soluble pathological tau, $n_x$ denotes the partial derivative of $n$ with respect to $x$, and $v(m,n)n$ is an active transport term. 
This velocity $v(m,n)$ is given by:
%
% \begin{equation}\label{velocity with exponentials}
% v(m,n)=v_ae^{\delta n}e^{-\varepsilon m}-v_r,
% \end{equation}
\begin{equation}\label{velocity}
v(m,n)=v_a(1+\delta n)(1-\varepsilon m)-v_r. 
\end{equation}
Here $v_a,v_r>0$ are the baseline anterograde and retrograde velocities of tau, respectively; $\delta$ is a nonnegative parameter governing the enhancement of kinesin processivity in response to soluble pathological tau, and $\varepsilon$ is a nonnegative parameter governing the reduction of kinesin processivity in response to insoluble pathological tau. 
The velocity term in \eqref{velocity} takes into account the fact that the propensity of soluble tau to travel in the anterograde direction is increased by a factor proportional to its concentration and decreased by a factor proportional to insoluble tau concentration. 

Because $m_{ij}$ represents the concentration of insoluble aggregates of tau, the diffusive and active transport flux of this species is defined to be 0. 
However, the interconversion between $m_{ij}$ and $n_{ij}$ on the edge through aggregation and fragmentation are given by:
\begin{equation}\label{equ edge}
\begin{aligned}
&\quad \begin{cases} 
\phi (m_{ij})_t=-\Gamma(m_{ij},n_{ij})&\text{in }(0,L)\setminus (x_3,x_4)\\
\Gamma(m_{ij},n_{ij}) = \beta m_{ij} - \gamma_{1} n_{ij}^{2} - \gamma_{2} n_{ij} m_{ij} \\
m_{ij}=0 &\text{in }(x_3,x_4)
\end{cases}\\
\phi (n_{ij})_t=&
\begin{cases} 
D(n_{ij})_{xx}\!+\!\Gamma(m_{ij},n_{ij})&\text{in }(0,x_1)\\
\lambda_1D(n_{ij})_{xx}\!+\!\Gamma(m_{ij},n_{ij})&\text{in }(x_1,x_2)\\
f D(n_{ij})_{xx}\!-\! (1\!-\!f) (v(m_{ij},n_{ij})n_{ij})_x \!+\!\Gamma(m_{ij},n_{ij})&\text{in }(x_2,x_3)\\
\lambda_2D(n_{ij})_{xx}&\text{in }(x_3,x_4)\\
D(n_{ij})_{xx}\!+\!\Gamma(m_{ij},n_{ij})&\text{in }(x_4,L),\\
\end{cases}
\end{aligned}
\end{equation}
where $f$ is the average fraction of soluble pathological tau that is undergoing diffusion as opposed to active transport at any given time \cite{Konzack2007,Cuchillo-Ibanez2008}, $\beta$ is the unimolecular rate of fragmentation, $\gamma_{1}$ is the bimolecular rate of soluble-soluble tau aggregation, $\gamma_{2}$ is the bimolecular rate of soluble-insoluble tau aggregation, $\lambda_{1}<1$ is the reduction of diffusivity in the AIS, and $\lambda_{2}<1$ is the reduction of diffusivity in the SC. In addition, $\phi>0$ is a small constant which represents the proportion between the slow and %$t$-scale and the 
fast time scales.  
Note that here we have generalized the Torok, \textit{et al.} model to allow for different aggregation rates for soluble-soluble and soluble-insoluble interactions as well as different diffusivities in the AIS and SC.

\begin{figure}
\centering
\includegraphics[width=1 \linewidth] {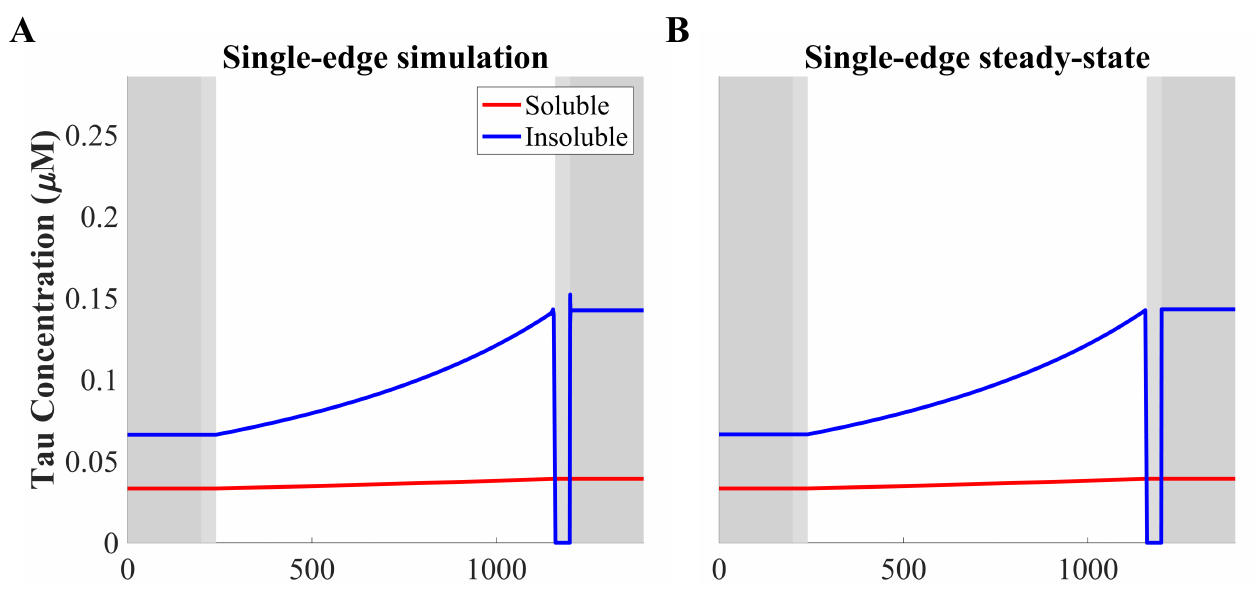}
\caption{Comparison of end-timepoint simulations versus the analytical steady-state for the single-edge model. A. Spatial distributions of soluble ($n$) and insoluble ($m$) tau across compartments after simulating the Torok \textit{et al.} single-edge mode to $t > 12$ months. B. Spatial distributions of soluble ($n$) and insoluble ($m$) tau across compartments for the derived steady-state solution of the Torok \textit{et al.} model.}
\label{fig:steadystate}
\end{figure}

In the original axonal transport model, the system was assumed to be closed at $x=0$ and $x=L$, and therefore Neumann zero-flux boundary conditions were imposed. 
Although we assume the overall network to be closed for the NTM, we require the transfer of mass between edges and nodes, and therefore these boundary conditions are insufficient.
Therefore, equations on the edges are completed by specifying biologically plausible initial conditions for $m_{ij}$ and $n_{ij}$ and the following \textit{Dirichlet} boundary conditions for soluble tau for all $t>0$:
\begin{equation}\label{BC edge}
n_{ij}(0,t)=N_i(t),\qquad n_{ij}(L,t)=N_j(t).
\end{equation}
 where $N_i(t)$ denotes the density (mass per unit volume) at vertex $P_i$ of soluble pathological tau protein.
We define the fluxes at the neuron-edge boundaries for 
node $P_{i}$ to be
\begin{align}\label{fluxes}
 J_{ij}^\phi(i,t)=-D (n_{ij})_x(0,t),\qquad J_{ji}^\phi(i,t)=-D (n_{ji})_x(L,t)
\end{align}
As we describe in more detail below, we also assume that all the biophysical processes on the edge occur on a ``fast'' time scale. 
Therefore, on the ``slow'' time scale of the network dynamics, $m_{ij}$ and $n_{ij}$ reach their steady-state distributions within the edge, which can be described by the following equations:
%\begin{equation}\label{equ edge, vanishing rho}
%\begin{cases} 
%m_{ij}=g(n_{ij})=\frac {\gamma_1 n_{ij}^2}{\beta-\gamma_2 n_{ij}}
%&\text{in }(0,L)\setminus (x_3,x_4)\\
%m_{ij}=0 &\text{in }(x_3,x_4)\\
%D(n_{ij})_{xx}=0&\text{in }(0,x_1)\\
%\lambda_1D(n_{ij})_{xx}=0&\text{in }(x_1,x_2)\\
%f D(n_{ij})_{xx}\!-\! (1\!-\!f) (v(g(n_{ij}),n_{ij})n_{ij})_x =0
%&\text{in }(x_2,x_3)\\
%\lambda_2D(n_{ij})_{xx}=0&\text{in }(x_3,x_4)\\
%D(n_{ij})_{xx}=0&\text{in }(x_4,L).\\
%\end{cases}
%\end{equation}
% coupled with boundary conditions as in \eqref{BC edge}. 
 
%For later purpose, we summarize the equation for $n_{ij}$ on the edge $e_{ij}$ as follows (see equations \eqref{BC edge} and \eqref{equ edge, vanishing rho}):
\begin{equation}\label{equ n, vanishing rho}
\begin{cases}
m_{ij}=g(n_{ij})=\frac {\gamma_1 n_{ij}^2}{\beta-\gamma_2 n_{ij}}&\text{in }(0,L)\setminus (x_3,x_4)\\
m_{ij}=0 &\text{in }(x_3,x_4)\\
\left(a(x)(n_{ij})_x+h(x,n_{ij})\right)_x=0 &\text{in }(0,L)\\
n_{ij}(0,t)=N_i(t), \quad n_{ij}(L,t)=N_j(t),
\end{cases}
\end{equation}
where 
\begin{equation}
\nonumber 
%\begin{cases}
 h(x,n_{ij})=\begin{cases}
-(1-f)v(g(n_{ij}),n_{ij})n_{ij} &  x \in (x_2,x_3) \\
 0 & \text{  otherwise}
\end{cases}
\qquad \;\;
%\end{cases}
%\end{equation}
%\begin{equation}
%\begin{cases}
 a(x)=\begin{cases}  
 D &\text{if } x \in (0,x_1)\\  %\bigcup  (x_4,L)\\
 D \lambda_1 &\text{if } x \in (x_1, x_2)\\
 f D  &\text{if } x\in (x_2, x_3)\\
 D \lambda_2 & \text{if } x \in (x_3, x_4)\\
 D &\text{if }x \in (x_4,L).
\end{cases}
%\end{cases}
\end{equation}
We observe that the flux of $n_{ij}$ on the edge $e_{ij}$ only depends on time:
\begin{equation}\label{constant flux}
J_{ij}(t)=-\left(a(x)(n_{ij})_x+h(x,n_{ij})\right).
\end{equation}
 As shown in \textbf{Figure~\ref{fig:steadystate}}, the steady-state distributions of $m_{ij}$ and $n_{ij}$ given by the equations above are equivalent to the those obtained by simulating the Torok, \textit{et al.} model for sufficiently long time (subject to the original Neumann zero-flux boundary conditions).

\subsection{Network Transport Model (NTM)}
%
%\mbcomment{Rewritten, referring to the single edge problem as a ``building block''.} 
% \jtcomment{What is this small parameter $\rho$? It also appears above in Equation 3 and it's unclear to me what it represents. Is it there to provide a theoretical upper bound on the amount of tau in the system so that the steady-state exists?}
% \vtcomment{We now have already defined $\rho$ in the single edge model }
Let the small parameter $\phi>0$ be fixed and let $M_i(t)$ and $N_i(t)$ denote the densities per unit volume at vertex $P_i$ of, respectively,  insoluble and soluble pathological tau protein, at time $t$. The equations for $M_i$ and $N_i$ are 
\begin{equation}\label{vertices, rho positive bis}
\begin{cases}
\phi \,M'_{i}=-\Gamma(M_i,N_i)\\
\phi \,N'_{i}=\frac 1{V_i}\underbrace{\sum_{j\ne i} 
\left(-c_{ij}J^{\phi}_{ij}(i,t) +c_{ji}J^\phi_{ji}(i,t)\right)}_{\text{incoming mass flow at compartment } P_i}+\Gamma(M_i,N_i),
\end{cases}
\end{equation}
%\mbcomment{Should the incoming flux contain the weights??$-c_{ij}J^\rho_{ij}(i,t) +c_{ij}J^\rho_{ji}(i,t)$ 
where $N_i'$ and $M_i'$ denote derivatives with respect to $t$, the reaction term $\Gamma(M_i,N_i)$ is defined by \eqref{equ edge}, $V_i$ is the volume of the brain compartment $P_i$, and $J^\phi_{ji}(i,t)$ and $J^\phi_{ji}(i,t)$ are the contributions to the incoming flux at $P_i$ from, respectively, single neurons of the edges $e_{ij}$ and $e_{ji}$, defined in \eqref{fluxes}.  We observe that, defining the total mass on the edge $e_{ij}$ as $c_{ij} \int_0^L (n_{ij} +m_{ij})(x,t) \;dx$ we assume the factor  $c_{ij}$
to be exactly the same as multiplying the fluxes $J_{ij}^\phi$ by $c_{ij}$.
In addition, due to the definition of mass flow, it is natural to choose the weights $c_{ij}$ to be the connectivity densities from the mouse structural connectome connectome (see \eqref{weights}) as these are constants that are proportional to the cross-sectional area of the white matter tracts represented by edges $e_{ij}$.  
%{\color{brown} Motivating our choice for $\rho N_i'$ in \eqref{vertices, rho positive bis}}

%{\color{brown} we assume that the weights $c_{ij}$ (see \eqref{weights}) are proportional to the cross-sectional area of the white matter tracts represented by edges $e_{ij}$. 
%Therefore, the incoming flux at $P_{i}$ is a summation of single-axon {\color{blue} edge?} fluxes $J_{ij}^\rho$ and $J_{ji}^\rho$ (see \eqref{fluxes}) weighted by $c_{ij}$ and $c_{ji}$, respectively.}
%\vtcomment{I suggest to skip this sentence.}
So for $\phi>0$ the NTM is described by system \eqref{vertices, rho positive bis}, completed by initial data for $N_i(0)$ and $M_i(0)$, and the single edge problem discussed in Section \ref{S edges}.

\subsection{Quasi-static approximation}\label{quasi-static}

%\mbcomment{Rewritten}

Since $\phi$ is a very small number we consider the limit $\phi\to 0$. 
This leads to a highly nontrivial singular perturbation problem and a rigorous mathematical treatment is far beyond the scope of the present paper. Instead we proceed formally and formulate a limit problem which can be considered as a quasi-static approximation of the model presented in the previous subsection.

%\subsubsection{Equations at vertices and edges}

Setting $\phi=0$ in the first equation of 
\eqref{vertices, rho positive bis}, 
we find that $M_i(t)$ is determined by $N_i(t)$:  
\begin{equation}\label{g(n)}
\Gamma(M_i,N_i)=0 \ \Leftrightarrow\ M_i=g(N _i):=\frac {\gamma_1 N_i^2}{\beta-\gamma_2 N_i}.
\end{equation}
Therefore also in the second equation of 
\eqref{vertices, rho positive bis}, the reaction $\Gamma$ 
formally disappears  
in the limit $\phi\to 0$, whence the second equation only should
account on the mass balance between the brain compartment $P_i$ and the 
edges $e_{ij}$ and $e_{ji}$. 
This mass balance contains two 
contributions. The first one are the constant mass fluxes entering $P_i$, 
i.e.\ $-J_{ij}(t)$ on $e_{ij}$ and $J_{ji}(t)$ on $e_{ji}$ (see \eqref{constant flux}.
%\vtcomment{I suggest to remove the weight from the flux }
The second contribution is caused by a feedback mechanism: 
a change of the Dirichlet condition $N_i(t)$ for the density 
$n_{ij}$ on $e_{ij}$ given by \eqref{BC edge},\eqref{equ n, vanishing rho}
%\mbcomment{add reference for $n_{ij}$ to single edge Section}
causes 
a change of the total mass on the edge which must be compensated by 
a change of the total mass at $P_i$. 
This leads to the following mass balance at $P_i$ at time $t$:
\begin{equation}\label{new mass balance}
\underbrace{V_i (N_{i}'+M_{i}')}_{\text{mass increase at $P_i$}}=\,
\underbrace{\sum_j \left(c_{ji}J_{ji}(t)-c_{ij}J_{ij}(t)\right)}_{\text{incoming mass flow at $P_i$}}\,
-\,\underbrace{\sum_j (C^i_{ij}(t)+C^i_{ji}(t)) N_i'(t)}_{\text{feedback mechanism}}.
%\gamma_{ji}R_{ji}^iN_{i}'-\gamma_{ij}S_{ij}^iN_{i}'.
\end{equation}
So we need to calibrate the coefficients $C^i_{ij}(t)$ and $C^i_{ji}(t)$ in order to obtain the correct mass balance at $P_i$.
%guarantee mass conservation.

Setting  $q_{ij}=\dfrac{\partial n_{ij}}{\partial t}$, the rate of change of the total mass on $e_{ij}$ is given by
\begin{equation}\label{rate mass chsange}
\begin{aligned}
&c_{ij}\int_0^L \frac{\partial}{\partial t}(n_{ij}+m_{ij})(x,t)\,dx
=c_{ij}\left(\int_0^L q_{ij}(x,t)\,dx+\int_{(0,x_3)\cup(x_4,L)}\frac {\gamma_1 n_{ij}(2\beta-\gamma_2  n_{ij})}{(\beta-\gamma_2  n_{ij})^2}q_{ij}(x,t)\,dx\right).
%\gamma_{ij}\int ( n_{ij})_t(x,t)\,dx,
\end{aligned}
\end{equation}
%where $\kappa_{ij}$ is a constant (for example, the average cross-section
%of the neuron bundle represented by the edge; a possible choice is 
%$\kappa_{ij}=c_{ij}$).
 Recalling  that $n_{ij}(t)$ satisfies \eqref{equ n, vanishing rho}, 
%the Dirichlet problem  (see \eqref{equ edge, vanishing rho})  %\mbcomment{check notation in single edge section} 
%$$
%\begin{cases}
%\left(a(x)(n_{ij})_x+h(x,n_{ij})\right)_x=0 \quad \text{in }(0,L)\\
%n_{ij}(0,t)=N_i(t), \quad n_{ij}(L,t)=N_j(t),
%\end{cases}
%$$
%{\color{blue}
%where 
%\begin{equation}
%\nonumber 
%%\begin{cases}
% h(x,n_{ij})=\begin{cases}
%-(1-f)v(g(n_{ij}),n_{ij})n_{ij} &  x \in (x_2,x_3) \\
% 0 & \text{  otherwise}
%\end{cases}
%\qquad \;\;
%%\end{cases}
%%\end{equation}
%%\begin{equation}
%%\begin{cases}
% a(x)=\begin{cases}  
% D &\text{if } x \in (0,x_1)\\  %\bigcup  (x_4,L)\\
% D \lambda_1 &\text{if } x \in (x_1, x_2)\\
% f D  &\text{if } x\in (x_2, x_3)\\
% D \lambda_2 & \text{if } x \in (x_3, x_4)\\
% D &\text{if }x \in (x_4,L),
%\end{cases}
%%\end{cases}
%\end{equation}
%}
we obtain that $q_{ij}(t)$ satisfies the linearized problem 
\begin{equation}\label{q_new}
\begin{cases}
\left(a(x)(q_{ij})_x+\frac{\partial h(x,n_{ij})}{\partial n}q_{ij}\right)_x=0 \quad \text{in }(0,L)\\
q_{ij}(0,t)=N_i'(t), \quad q_{ij}(L,t)=N_j'(t).
\end{cases}
\end{equation}
The linearity of this problem makes it possible to distinguish the contributions to the rate of total mass change on $e_{ij}$ caused by 
$N'_i$ and $N_j'$: 
\begin{equation}\label{linear combination}
q_{ij}(x,t)=N'_i(t)q_{ij}^i(x,t)+N'_j(t)q_{ij}^j(x,t),
\end{equation}
where, for all fixed $t$, the functions $q_{ij}^i(t)$ and $q_{ij}^j(t)$ satisfy, respectively,
\begin{equation}\label{q_new1}
\begin{cases}
\left(a(x)(q_{ij}^i)_x+\frac{\partial h(x,n_{ij})}{\partial n}q^i_{ij}\right)_x=0 \quad \text{in }(0,L)\\
q^i_{ij}(0,t)=1, \quad q^i_{ij}(L,t)=0
\end{cases}
\qquad
\begin{cases}
\left(a(x)(q_{ij}^j)_x+\frac{\partial h(x,n_{ij})}{\partial n}q^j_{ij}\right)_x=0 \quad \text{in }(0,L)\\
q^j_{ij}(0,t)=0, \quad q^j_{ij}(L,t)=1.
\end{cases}
\end{equation}
Therefore \eqref{linear combination} strongly suggests to define, in \eqref{new mass balance},
\begin{equation}\label{Cij}
C_{ij}^i(t)=c_{ij}\left(\int_0^L q_{ij}^i(x,t)\,dx+\int_{(0,x_3)\cup(x_4,L)}\frac {\gamma_1 n_{ij}(2\beta-\gamma_2  n_{ij})}{(\beta-\gamma_2  n_{ij})^2}q_{ij}^i(x,t)\,dx\right).
\end{equation}
Similarly we set
\begin{equation}\label{Cji}
C_{ji}^i(t)=c_{ji}\left(\int_0^L q_{ji}^i(x,t)\,dx+\int_{(0,x_3)\cup(x_4,L)}\frac {\gamma_1 n_{ji}(2\beta-\gamma_2  n_{ji})}{(\beta-\gamma_2  n_{ji})^2}q_{ji}^i(x,t)\,dx\right).
\end{equation}

\subsubsection{The complete quasi-static NTM}\label{complete model}

%\mbcomment{A first attempt}

We briefly reassume the mathematical network-transport model under the quasi-static assumption. 

Given a directed graph with vertices $P_i$ (the brain compartments with volumes $V_i$) and edges $e_{ij}$ 
with weights $c_{ij}$ (bundles of connecting neurons between compartments $P_i$ and ($P_j$), 
the quasi-static NTM consists of equations for the densities $N_i(t)$ and $n_{ij}(x,t)$ for soluble Tau at $P_i$ and, respectively, on $e_{ij}$. 
The densities $M_i$ and $m_{ij}$ for insoluble Tau then follow from the quasi-static assumption: $M_i=g(N_i)$ and $m_{ij}=g(n_{ij})$ (see 
\eqref{vertices, rho positive bis}).

We start from given initial data $N_{0i}$ for $N_i$ at all brain compartments $P_i$:
$$
N_i(0)=N_{0i} \in [0, \tfrac \beta{\gamma_2})
$$
(in all numerical simulations we always set $\gamma_2=0$, so we simply have $N_{0i}\ge 0$). 
The model has two major ingredients:
\begin{itemize}
\item given the concentrations $N_i(t)$ at $P_i$ at a fixed time $t\ge 0$, the concentrations $n_{ij}(x,t)$ on the edges are determined by the 
single-edge model \ref{S edges};
%\mbcomment{add reference to single edge Section}
the mass flux $c_{ij}J_{ij}(t)$ on $e_{ij}$ from $P_i$ to $P_j$ 
is constant along the edge and is completely determined by $n_{ij}(x,t)$;
\item $N_i$ is determined by the differential equation \eqref{new mass balance} with initial condition $N_{i0}$, 
where $M_i=g(N_i)$ and where the coefficients $C_{ij}^i(t)$ and $C_{ij}^i(t)$ are given by 
\eqref{Cij} and \eqref{Cji}:
\begin{equation}\label{quasi-static model}
\left(V_i \left(1+\frac {\gamma_1 N_{i}(2\beta-\gamma_2 N_{i})}{(\beta-\gamma_2 N_{i})^2}\right)+
\sum_j \left(C_{ij}^i+C_{ji}^i\right)\right)N_{i}'=\sum_j \left(c_{ji}J_{ji}-c_{ij}J_{ij}\right).
\end{equation}
\end{itemize}

If $\gamma_2=0$ (the choice in all numerical simulations presented below), the NTM quasi-static problem turns out to 
possess a well-defined solution which is uniquely determined by its initial data $N_{0i}\ge 0$. The proof goes beyond the scope of the present paper and will
be presented in a future paper. % ([Emilia e.a.]). 
There we shall also treat the more general case of positive $\gamma_2$, which requires a more refined
analysis to handle the singularity $n=\beta/\gamma_2$ in the relation $m=g(n)$.

Finally we mention that in Torok, \textit{et al.}\cite{Torok2021}, computational tests on a 
single edge, with no flux instead of Dirichlet boundary conditions,
strongly suggest that for positive $\phi$ the solutions of the PDE \eqref{equ edge}, 
rapidly converge to equilibrium for $t>0$ if $\phi$ is very small. 
This provides a first computational justification for the introduction of 
the quasi-static approximation.

$$
%M_i(0)=M_{0i}=\frac {\gamma_1 N_{0i}^2}{\beta-\gamma_2 N_{0i}}.
$$
\section{Numerical simulations}
%
%
%
%
%

% \subsection{Software and numerical implementation}\label{software}
\subsection{Model implementation}

We implemented the quasi-static approximation to the NTM (\eqref{Cij}-\eqref{quasi-static model}) in MATLAB version 2022b. 
%\jtcomment{get correct eq. refs.}
To make this problem more computationally feasible, the equation in \eqref{quasi-static model} has been discretized using a first-order Euler method and the integrals in \eqref{Cij},\eqref{Cji} have been approximated through a  trapezoidal method (see \textbf{Appendix A}).
The steady-state solution $n_{ij}, i \sim j$ on each neuronal bundle has been calculated by means of the MATLAB solver \texttt{ode45} which is based on an explicit Runge-Kutta method of order (4,5).  
The mass fluxes $J_{ij}$ and $J_{ji}$ as well as the integration constants in \eqref{q_new} have been calculated at each model time on each neuronal bundle through  shooting procedures relying on the use of MATLAB's nonlinear solver \texttt{fsolve}. 
We ran simulations in parallel using the computational resources provided by the University of California, San Francisco. 
\subsection{Computational experiments}
In order to provide a global picture of the evolution of tau pathology on the hippocampal subnetwork, we first define the total tau burden at each the vertex of the network as:
\begin{align}\label{tot_tau}
\nonumber &\tau(t):=(\tau_i(t),\dots,\tau_h(t)) \quad \text{and} \\ 
&\tau_i(t) = N_i(t)+M_i(t) \;\; i=1,\dots,h
\end{align}
where $h$ is the number of vertices of the graph $G$ (see \textbf{Section 2.1}).
We chose the lateral part of the left entorhinal cortex (ECl\textsuperscript{L}) as our initiation site for tau pathology, as it has been previously identified as one of the earliest and most affected regions in AD \cite{Braak1991}.
Mathematically, this can be expressed in the following way:
\begin{align}\label{in_data}
\nonumber &\tau_i(0)=0 \quad i \neq i_{seed}\\
&\tau_i(0)=0.02 \quad i= i_{seed}
\end{align}
where $i_{seed}$ is the index corresponding to %the 
the (ECl\textsuperscript{L}). 
Without loss of generality, we set the initial tau burden to be 0.02 at the seeding location.

To explore the model dynamics, we ran a grid search on five key model parameters: $\lambda_1$, $\lambda_2$, $\gamma_1$, $\delta$, and $\epsilon$. 
It has been previously shown that the ratio of aggregation rate to fragmentation rate governs the distribution of tau for the single-edge model \cite{Torok2021}, and therefore we hold the fragmentation rate, $\beta$, constant. 
As mentioned above in \textbf{Section \ref{quasi-static}}, the system exhibits a singularity for $\gamma_2 > 0$; therefore, we set $\gamma_2 = 0$ for all simulations. 
For simplicity, we also only consider the case where $\lambda_1 = \lambda_2$.

We consider the dynamics of the NTM in three important cases: 1) varying the strength of the diffusion barrier separating vertices and edges (governed by $\lambda_1$ and $\lambda_2$); 2) varying the degree of directional bias (governed by $\delta$ and $\epsilon$); and 3) varying the aggregation rate (governed by $\gamma_1$). 
For each experiment, we varied the parameter(s) of interest while holding all other parameters constant.
In particular, we sought to examine the influence of these parameters on two key aspects of tau dynamics: 1) the global rate of spread tau on the network, and 2) the staging of regional tau pathology, which can be described in terms of the ordering of regions by their peak tau concentration and the amount of time it takes for each region to reach those peak concentrations (\textit{arrival times}). 
\subsubsection{Varying diffusion barrier strength}
\begin{figure}
\centering
\includegraphics[width=1\linewidth]{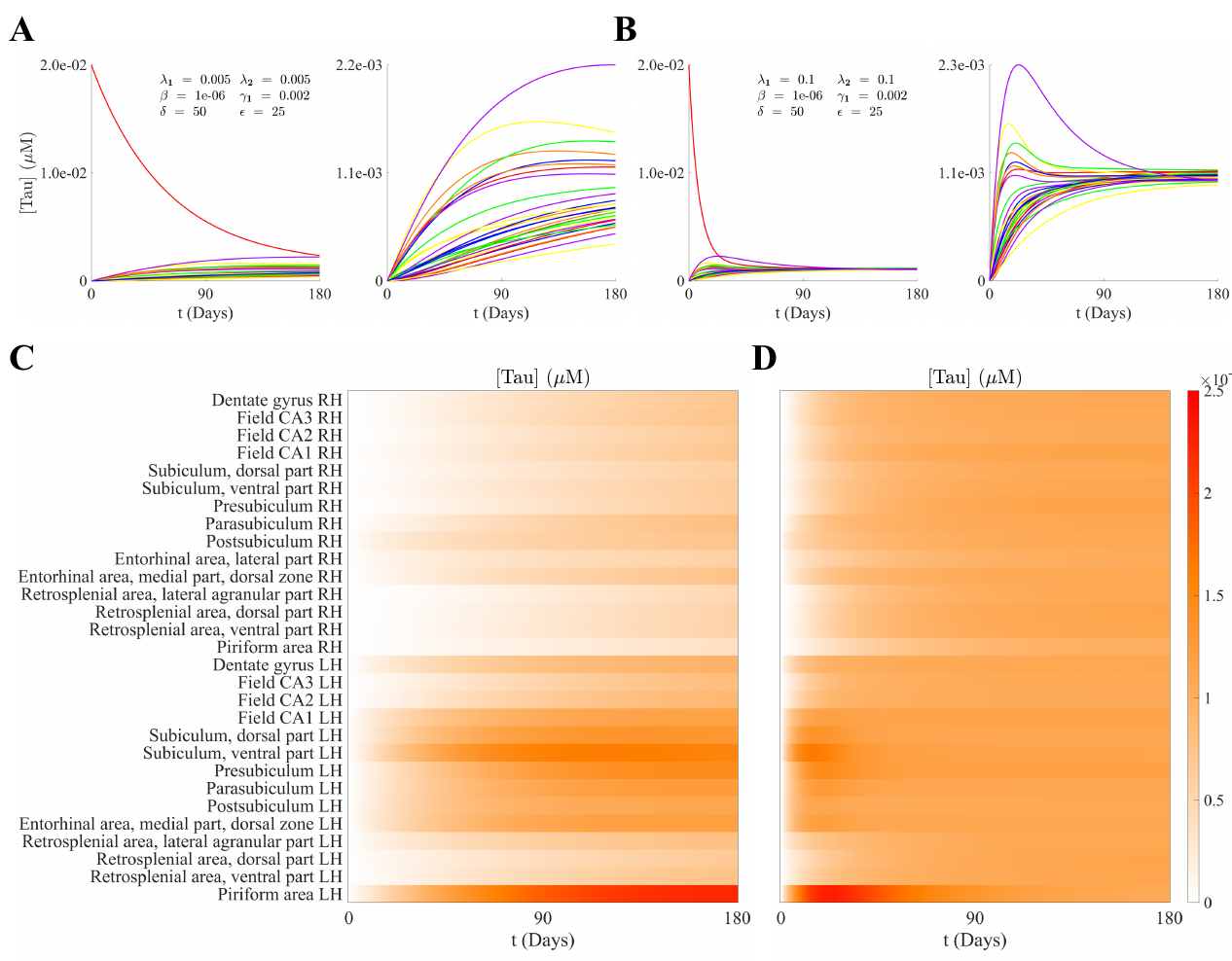}
\caption{Effect of $\lambda$ on NTM simulations. A. Total concentration of tau plotting against time in each of the 30 regions for the low lambda condition ($\lambda_{1} = \lambda_{2} = 0.005$). Left panel includes the seed region (left lateral part of the entorhinal cortex), right panel excludes the seed region. B. Tau concentration over time for the high lambda condition ($\lambda_{1} = \lambda_{2} = 0.1$). C. Heatmap representation of the per-region simulations shown in A., with the seed region excluded. D. Heatmap representation of the per-region simulations shown in B., with the seed region excluded.
}
\label{fig:lambdamodels}
\end{figure}

\textbf{Figure~\ref{fig:lambdamodels}} shows the spatiotemporal evolution of tau in each of the 30 regions of the hippocampal subnetwork for two sets of values for $\lambda_1$ and $\lambda_2$.
For $\lambda_1 = \lambda_2 = 0.005$, spread along the network is relatively slow (\textbf{Figure~\ref{fig:lambdamodels}A},\textbf{C}), where the tau burden in all unseeded seeds grows roughly monotonically with no signs of plateauing over the 180-day simulation. 
However, by increasing $\lambda_1$ and $\lambda_2$ to 0.1, we observe tau burden rapidly converges to nearly equivalent concentrations in all regions in this same time window, with clear peaks in the more affected unseeded regions (\textbf{Figure~\ref{fig:lambdamodels}B},\textbf{D}). 
Given that the $\lambda$ parameters exert their effects by reducing the effective diffusivity of soluble tau between axonal and somatodendritic compartments (see \textbf{Section 2.1}), we expected that reducing the values of $\lambda_1$ and $\lambda_2$ should have slowed the overall rate of tau spreading on the network, as the NTM demonstrated.
Moreover, we found that the $\lambda$ parameters do not appear to affect the order of arrival times of tau in each region. 
Comparing between \textbf{Figure~\ref{fig:lambdamodels}C} and \textbf{D}, while the rate of pathology spread is slower for the smaller value of $\lambda$, it otherwise very closely mirrors the spread of tau for the higher $\lambda$. 
We conclude that $\lambda$ strongly influences the global rate of pathology spread (providing a ``bottleneck''-like effect), but not the staging of pathology.

\subsubsection{Varying molecular motor rates}
As mentioned in \textbf{Section ~\ref{S edges}} above, the parameters $\delta$ and $\epsilon$ affect the processivity of kinesin along axonal microtubules as a function of pathological tau species, with soluble tau \textit{increasing} kinesin processivity through $\delta$ and insoluble tau \textit{decreasing} kinesin processivity through $\epsilon$. 
In terms of tau migration within a given axon, the effect of changing $\delta$ and $\epsilon$ is to confer an anterograde or retrograde spread bias, respectively~\cite{Torok2021}.

\begin{figure}
\centering
\includegraphics[width=1 \linewidth] {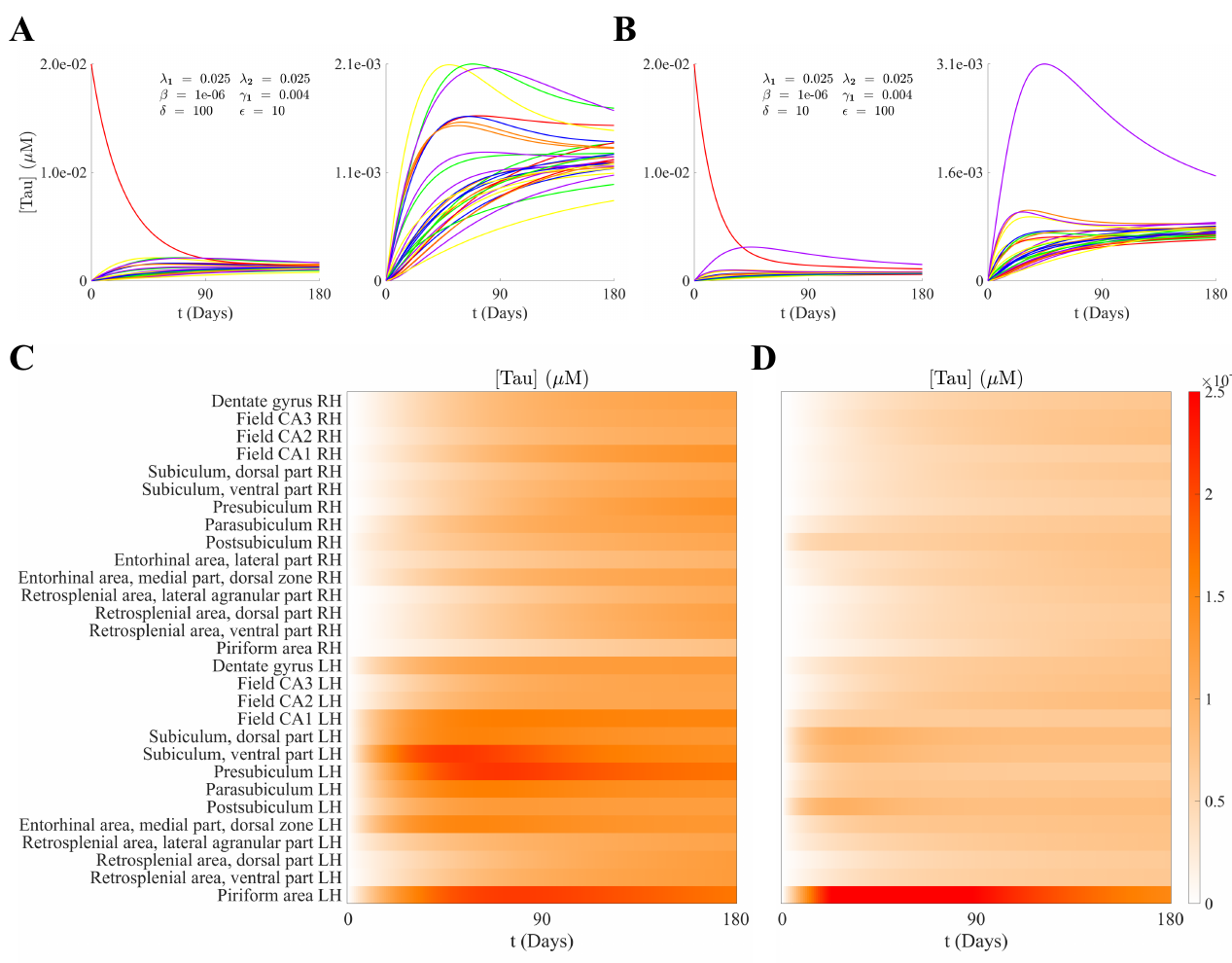}
\caption{Effect of $\delta$ and $\epsilon$ on NTM simulations. A. Total concentration of tau plotting against time in each of the 30 regions for the anterograde-biased condition ($\delta = 100, \epsilon = 10$). Left panel includes the seed region (left lateral part of the entorhinal cortex), right panel excludes the seed region. B. Tau concentration over time for the retrograde-biased condition ($\delta = 10, \epsilon = 100$). C. Heatmap representation of the per-region simulations shown in A., with the seed region excluded. D. Heatmap representation of the per-region simulations shown in B., with the seed region excluded.}
\label{fig:antretmodels}
\end{figure}

\begin{figure}
\centering
\includegraphics[width=1 \linewidth] {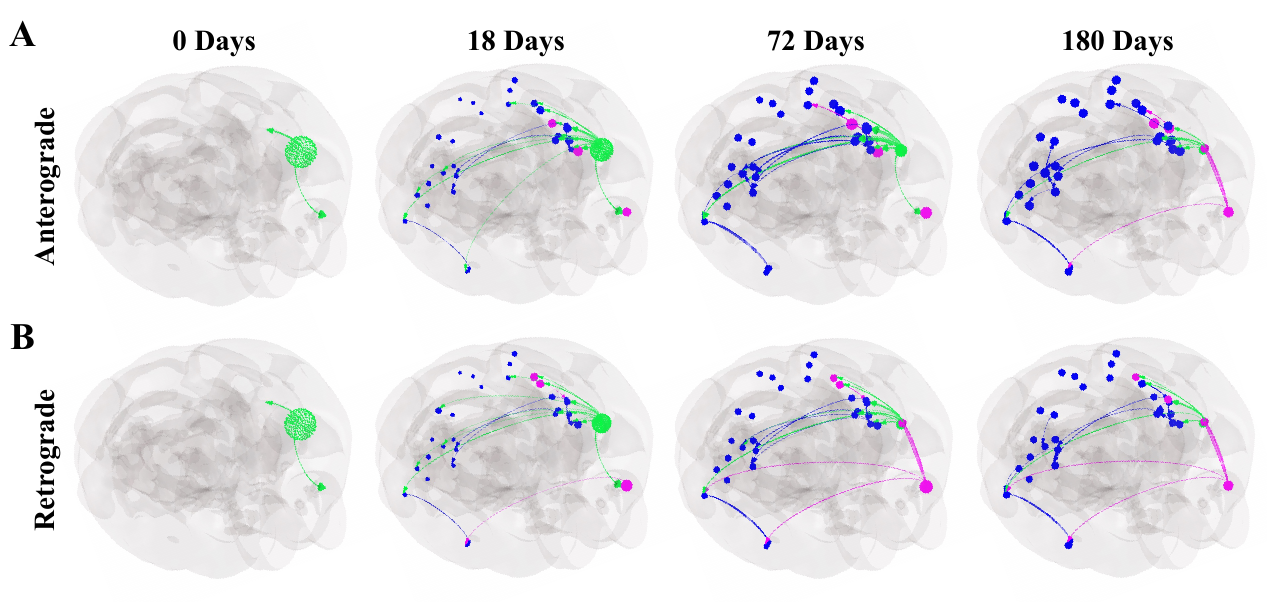}
\caption{Three-dimensional visualization of anterograde and retrograde biases A. Total concentration of tau over time, visualized as spheres centered in each region for the anterograde simulation shown in Figure~\ref{fig:antretmodels}A, where sphere size is proportional to the amount of pathology. Arrows represent the directions and strengths of the upper 10\% of fluxes between regions at each time point. All regions are colored in blue, with the exception of the seed region (in green) and the three non-seed regions with the most tau pathology at each time point (in magenta). B. Total concentration of tau per region in the retrograde simulation (Figure~\ref{fig:antretmodels}B), with the same visualization conventions as in A.}
\label{fig:antretgb}
\end{figure}

In contrast with $\lambda$, changing $\delta$ and $\epsilon$ in the NTM dramatically changed the spreading of pathology on the hippocampal subnetwork (\textbf{Figure~\ref{fig:antretmodels}}). 
In the anterograde-biased case ($\delta = 100, \epsilon = 10$), we observed that the three regions most strongly affected by tau pathology were, in order of arrival time: the left ventral part of the subiculum, the left presubiculum, and the left piriform area (\textbf{Figure~\ref{fig:antretmodels}A},\textbf{C}).
However, the in the retrograde-biased simulation ($\delta = 10, \epsilon = 100$), tau spread predominantly to the left piriform area only (\textbf{Figure~\ref{fig:antretmodels}B},\textbf{D}). 
The global rate of spread, which can be observed most clearly by examining the decrease of tau concentration in the seeded region, were similar between the anterograde-biased and retrograde-biased experiments 
(\textbf{Figure~\ref{fig:antretmodels}A},\textbf{B}, left panels, red lines). 
We visualize these differences in the brain at selected time points in \textbf{Figure~\ref{fig:antretgb}}, where we represent the tau burden each region as spheres with radii proportional to tau concentration and the strength and direction of tau fluxes between regions with arrows. 
We observed that, while pathology and fluxes at $t = 0$ were very similar for the anterograde and retrograde cases (\textbf{Figure~\ref{fig:antretmodels}A},\textbf{B}, respectively), the evolution of tau pathology grew increasingly divergent between these simulations over time. 
Notably, there was little overlap between simulations in terms of the unseeded regions most strongly affected by tau at each time point (magenta spheres). 
We conclude that the $\delta$ and $\epsilon$ parameters largely exert their effects by changing which regions were most strongly affected by tau as well as their arrival times.

\begin{figure}
\centering
\includegraphics[width=1 \linewidth] {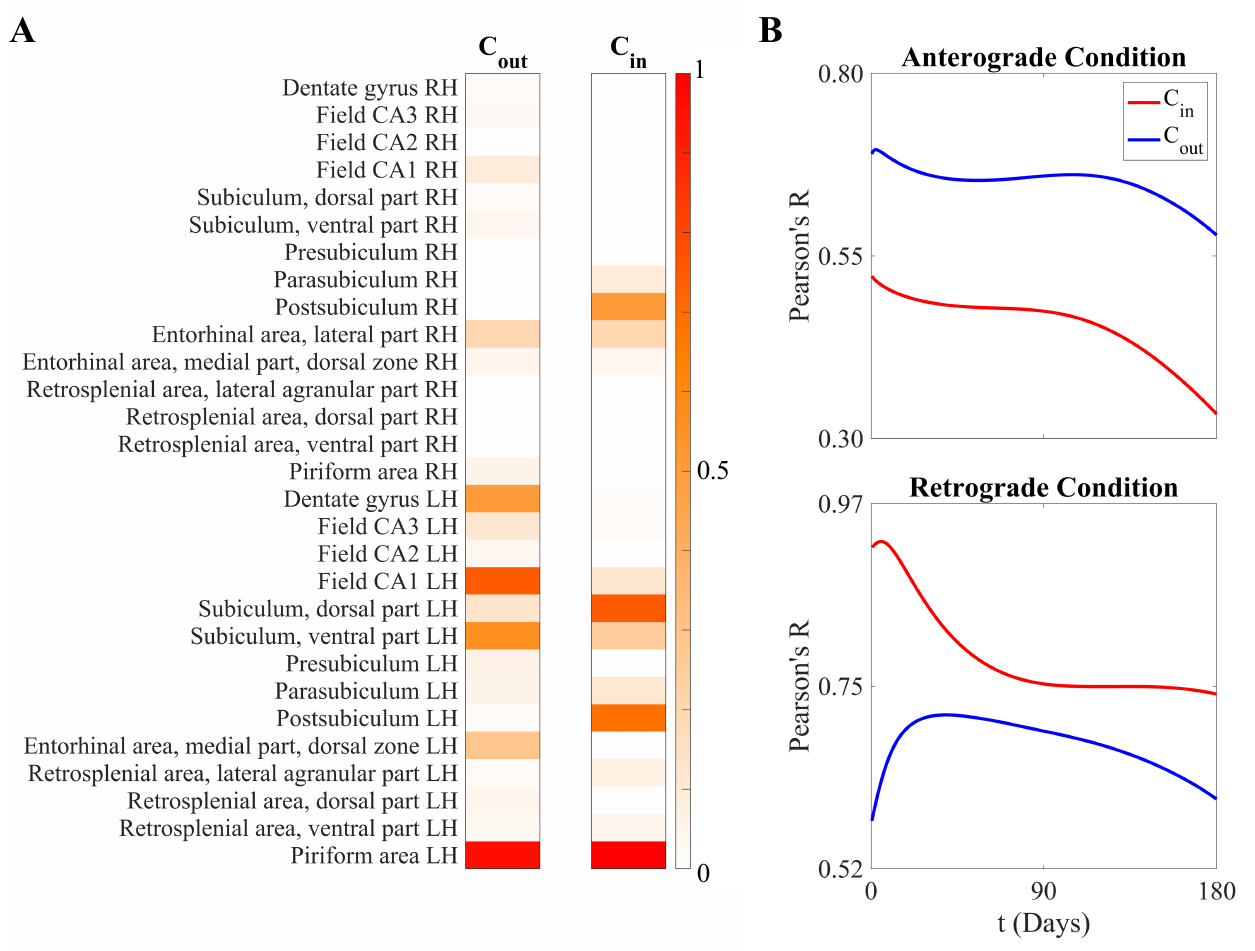}
\caption{Comparison of connectivity with respect to the seeded region and the anterograde/retrograde simulations. A. Outgoing ($\text{C}_{\text{out}}$) and incoming ($\text{C}_{\text{in}}$) connectivity with respect to the seeded region (left lateral entorhinal cortex). B. Spatial correlations with outgoing and incoming connectivity with respect to the seeded region and the spatiotemporal evolution of the anterograde-biased (top) and retrograde-biased (bottom) simulations. }
\label{fig:antretseedconn}
\end{figure}

Given the apparent impact of $\delta$ and $\epsilon$ on the staging of tau pathology, we sought to confirm that tau was spreading predominantly in the anterograde and retrograde directions for the corresponding values of these parameters.
In \textbf{Figure~\ref{fig:antretseedconn}A}, we depict the outgoing and incoming connectivity density with respect to the seed region, which we denote by $\text{C}_{\text{out}}$ and $\text{C}_{\text{in}}$, respectively. 
The left piriform area is the most strongly connected region to the left lateral part of the entorhinal cortex with respect to both outgoing and incoming connectivity, which explains why the left piriform area exhibits strong tau pathology for both anterograde-biased and retrograde-biased simulations. 
However, there are numerous differences elsewhere in the hippocampal subnetwork; for instance, the left dentate gyrus receives many projections from the ECl\textsuperscript{L} but largely does not project to it ($\text{C}_{\text{out}} > \text{C}_{\text{in}}$), whereas the left dorsal part of the subiculum projects strongly to the ECl\textsuperscript{L} but largely does not receive projections from it ($\text{C}_{\text{out}} < \text{C}_{\text{in}}$; \textbf{Figure~\ref{fig:antretseedconn}A}).
We hypothesized that the tau distributions would be more strongly associated with outgoing connectivity in the anterograde-biased case and more strongly associated with incoming connectivity in the retrograde-biased case.
Indeed, we observed that tau distributions in the anterograde-biased case were more correlated with ($\text{C}_{\text{out}}$) than ($\text{C}_{\text{in}}$), while the reverse was true in the retrograde-biased case (\textbf{Figure~\ref{fig:antretseedconn}B}). 
Therefore, we conclude that increasing or decreasing $\delta$ relative to $\epsilon$ causes tau to preferentially migrate in the anterograde or retrograde directions, respectively.

\subsubsection{Varying aggregation rate}
\begin{figure}
\centering
\includegraphics[width=1 \linewidth] {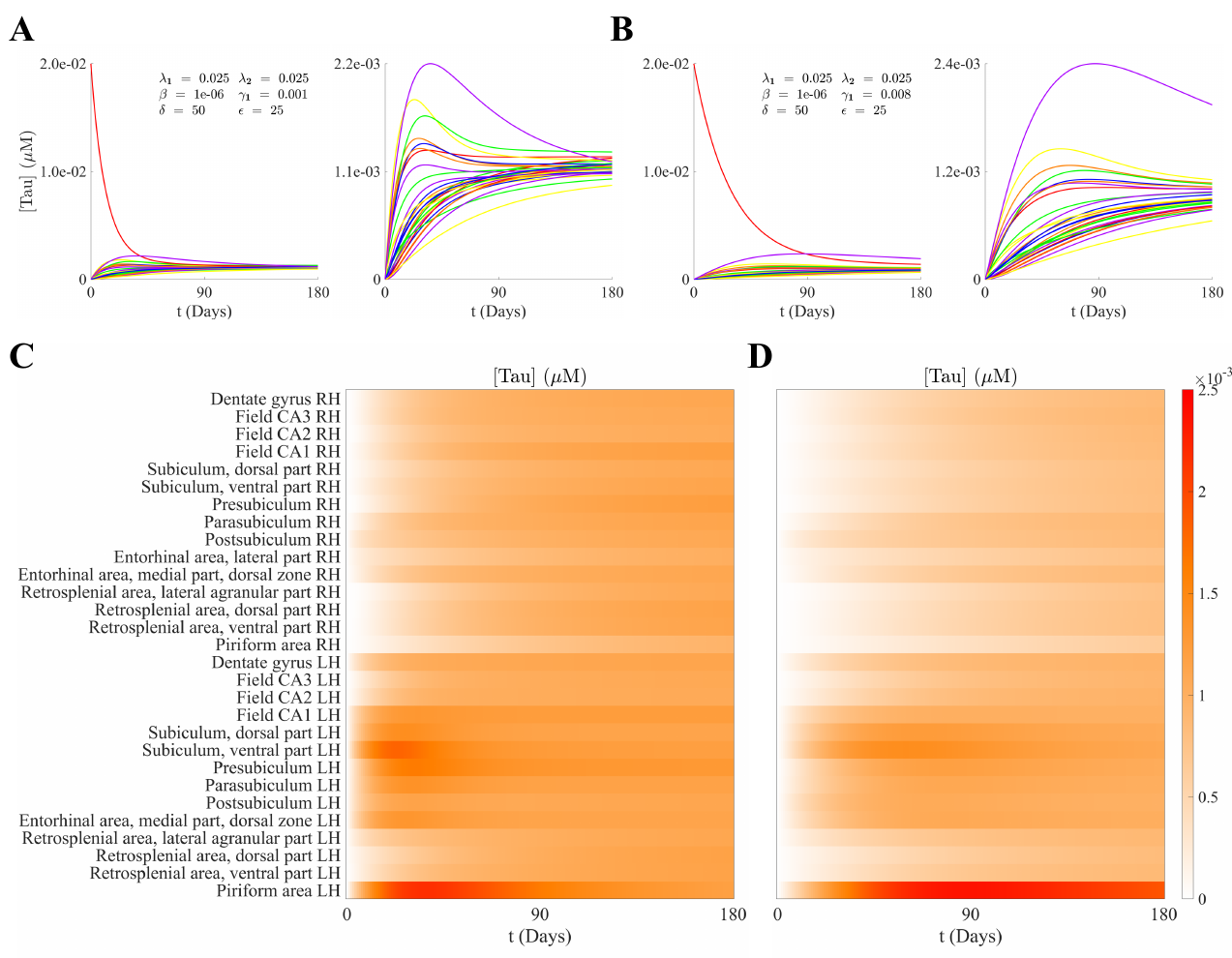}
\caption{Effect of $\gamma_{1}$ on NTM simulations. A. Total concentration of tau plotting against time in each of the 30 regions for the low $\gamma_{1}$ condition ($\gamma_{1} = 0.001$). Left panel includes the seed region (left lateral part of the entorhinal cortex), right panel excludes the seed region. B. Tau concentration over time for the high $\gamma_{1}$ condition ($\gamma_{1} = 0.008$). C. Heatmap representation of the per-region simulations shown in A., with the seed region excluded. D. Heatmap representation of the per-region simulations shown in B., with the seed region excluded.}
\label{fig:gamma1models}
\end{figure}

Lastly, we explored the effect of modulating the aggregation rate, $\gamma_1$, on the spread along the hippocampal subnetwork. 
Most notably, increasing $\gamma_1$ from 0.001 to 0.008 resulted in a reduction in the global rate of spread along the network (\textbf{Figure~\ref{fig:gamma1models}}), similar to the effect of decreasing $\lambda$.
Given that increasing $\gamma_1$ increases the amount of insoluble tau relative to soluble tau (see \textbf{Equation~\ref{equ n, vanishing rho}}) and both models were instantiated with the same amount of total tau, we attribute the slower rate of spread in the high $\gamma_1$ to an overall lower concentration of available soluble tau to migrate between regions. 
Interestingly, we observed that $\gamma_1$ also exerts a subtle effect on the staging of tau pathology. 
Although the left piriform area and the left ventral part of the subiculum were the unseeded regions exhibiting the highest tau burden for both values of $\gamma_1$, we observed that the left presubiculum (\textbf{Figure~\ref{fig:gamma1models}A}, green line) reached a higher peak tau concentration than the left dorsal part of the subiculum (\textbf{Figure~\ref{fig:gamma1models}A}, orange line) for $\gamma_1 = 0.001$, while this situation was reversed for $\gamma_1 = 0.008$.
Because the relative strengths of the kinesin processivity parameters $\delta$ and $\epsilon$ depend on the available amounts of soluble and insoluble tau, by increasing $\gamma_1$ and pushing the system towards insoluble tau, the staging of tau pathology on the network also changes.

% The effect of the aggregation parameter $\gamma_1$ on total Tau dynamics is tested for low value of $\gamma_1\; (0.001)$ versus high value of $\gamma_1 \; (0.008)$. (see Figure \ref{fig:aggmodels})
% %Longitudinal graphs of total Tau concentration as well as heat maps are plotted. 
% Differences in total Tau spatial distribution over time arise in the two cases. 
% For $\gamma_1=0.001$, the spread of pathological Tau from the seeding region affect mainly the Piriform area LH, Subiculum LH, Presubiculum LH at earlier times, while at later times the regions exhibiting a higher Tau load are Field CA1 RH, Presubiculum RH, Field CA1 RH, Presubiculum LH, although the distribution of Tau pathology results rather uniform in the whole network. 
% On the other hand, for $\gamma_1=0.008$, Tau pathology spreads mainly towards the Piriform area LH. 
% Tau distributions for low/high $\gamma_1$ resemble respectively those  of the anterograde/ retrograde dynamics. Indeed, the parameter $\gamma_1$  influence the equilibrium between soluble and insoluble Tau, and consequently modulate the effect of $\delta/ \varepsilon$ on the transport velocity $v$. Low value of $\gamma_1$ implies more soluble tau than insoluble tau in the system, then the effect of the parameter $\delta$ is increased pushing the dynamics more in the anterograde direction. Analogously, high values of $\gamma_1$ are linked with more insoluble than soluble tau, implying a leading impact of the parameter $\varepsilon$ on the transport velocity and consequently a governed-retrograde dynamics.
%
%
%
%
%
\section{Discussion}
%
% AR Added this
Mathematical modeling of tau propagation in Alzheimer's disease in a ``prion-like" fashion on the structural connectome of the brain is a recent yet powerful advance which can provide a platform for testing novel hypotheses regarding pathophysiology without the need for expensive and time-consuming experimental studies. 
Prior work has established the efficacy of modeling trans-synaptic spread of tau via either anisotropic diffusion in inhomogeneous media \cite{Kuhl2019} or diffusion on the connectivity graph \cite{Raj2012, Torok2018, Mezias2020a}. 
Many extensions and applications of these models are available ~\cite{Fornari2019,Iturria-Medina2014,Schafer2020,Weickenmeier2018,Weickenmeier2019,Bertsch2023}. 
However, more recent work has highlighted the important role of active axonal transport of tau via molecular motors attached to microtubules in either anterograde or retrograde direction, and the interactions between tau hyperphosphorylation, microtubules and the motor proteins ~\cite{Cuchillo-Ibanez2008,Rodriguez-Martin2013,Stern2017,Alonso1994,Alonso1997,Horvat2023}. 
Together, these effects lead to aberrant axonal transport and the mis-sorting of tau into the neuronal somatodendritic compartment~\cite{Ballatore2007,Zempel2017}. 
These transport-related effects were recently modeled mathematically for a single axon by Torok, \textit{et al.} who proposed a two-species, multicompartment model to explore the interactions between the pathological axonal transport of tau and the formation and breakdown of insoluble tau aggregates~\cite{Torok2021}. 
It was shown that active transport and its interaction with other protein kinetic processes is a key aspect of tau progression, which cannot be captured by prior diffusion models. 

Hence we undertook the current study to fill the gap between brain-wide diffusion on the connectivity graph and transport models which are localized to a single edge of that graph. 
We propose a macroscopic model called the Network Transport Model (NTM) which combines the dynamics of soluble and insoluble tau in the gray-matter regions and the Torok, \textit{et al.} axonal transport model~\cite{Torok2021} in the white matter tracts. 
The NTM was designed to simulate the dynamics of soluble and insoluble tau in terms of the diffusion-advection and aggregation-fragmentation processes as in ~\cite{Torok2021}, but at the network level. Although the full mathematical problem is 
computationally intractable, we developed a quasi-static approximation which maintains the basic properties of the full NTM by assuming a separation of time scales and solving the resulting mass balance problem.
% This required some key mathematical innovations. We separated the full network dynamics into two time scales: a fast time scale (order of hours to days) whereby transport processes in individual edges are established; followed by a slow time scale (order of days to months) over which the full network couples to these individual edge processes via slow exchange within network nodes. We imposed a a steady state approximation of the PDE solution at the fast scale at the terminals of the edge and solved it using tractable numerical integration. At the slow time scale we formalized the detailed mathematical description of the mass transfer between edges and nodes under a local mass balance problem. 

The resulting quasi-static NTM model was evaluated extensively using numerous simulations on the mouse mesoscale connectome obtained from viral tracing experiments \cite{Oh2014}. 
We performed three key experiments: 1) varying the strength of the diffusivity barrier, 2) changing the directionality bias, and 3) changing the aggregation rate.
We found that the system evolved much faster when the effective diffusivity within the single-edge model was higher (i.e., higher $\lambda_{1,2}$, but the staging of pathology did not change. 
By contrast, by changing the tau-kinesin parameters $\delta$ and $\epsilon$ predominantly affected the staging of pathology and demonstrably shifted the system towards either an anterograde or retrograde bias.
The aggregation rate showed a combination of effects, where higher aggregation rates were associated with a slower rate of spread and subtle differences in regions most strongly affected by tau.

Taken together, the above computational experiments demonstrate the rich dynamics of the NTM as a function of microscopic kinetic and transport parameters, and is able to sustain a far more diverse range of behaviors than was possible with previous diffusion-based models. 
Of particular import is the ability of the NTM to govern directionally biased flows on the connectome, an aspect that  has received almost no theoretical attention previously. 
This quasi-static NTM is therefore capable of delivering more robust insights into how specific tau species may differentially propagate in the brain, which has important clinical implications for our understanding and treatment of AD and other tauopathic dementias. 
Future modeling and empirical studies are required to test the model on empirical data, establish realistic parameter regimes and test specific hypotheses. 
Additional effort will also be needed to convert the current model into one that can be practically applied to real cases quickly and efficiently. 
%
% AR portion ends

\bigskip

 \noindent {\bf Acknowledgement.} The authors are strongly indebted to Prof.\ Bruno Franchi for stimulating their collaboration.\\
M.B.\ and V.T.\ acknowledge the MIUR Excellence Department  Project \texttt{MatMod@TOV} awarded 
to the Department of Mathematics, University of Rome Tor Vergata.

\bibliographystyle{abbrv}
\bibliography{Pathology_Master}
\clearpage
\section{Figures and Tables}
\begin{table}[h]
\begin{tabular}{ l l l } 
 \hline
 \bf{Symbol} & \bf{Description} & \bf{Remark} \\ 
 \hline
\\
 $x$ & Space  &  {\small $x \in [0, L]$, where $L$ is the total size of the system in $\mu$m}\\
%  &  & {\small is the length of the fiber.}\\ 
 $t$ & Time & {\small $t \in [0,T]$, where $T$ is the time to reach steady state in s} \\
%   &  & {\small is the time to reach steady state.}\\ 
$n_{ij}$ & Soluble tau & {\small Number of monomeric units per volume within} \\
   &  & {\small a single set of neurons within edge $e_{ij}$, in $\mu$M}\\ 
$m_{ij}$ & Insoluble tau & {\small Number of units of aggregates per volume within} \\
   &  & {\small a single set of neurons within edge $e_{ij}$, in $\mu$M}\\
$D$ & Theoretical diffusivity of $n$ & {\small Estimated to be 12 $\mu$m$^{2}$/s*} \\
$f$ & Diffusing fraction of $n$ & {\small Estimated to be 0.92*} \\
$v_{a}$ & Native ant. transport velocity of $n$ & {\small Estimated to be 0.7 $\mu$m/s*} \\
$v_{r}$ & Native ret. transport velocity of $n$ & {\small Estimated to be 0.7 $\mu$m/s*} \\
$\beta$ & Fragmentation rate of $m$ & {\small Unimolecular process by which $m \rightarrow n$}\\
$\gamma_1$ & Aggregation rate, $n+n$ & {\small Bimolecular process by which $n \rightarrow m$}\\
$\gamma_2$ & Aggregation rate, $n+m$ & {\small Bimolecular process by which $n \rightarrow m$}\\
$\delta$ & Ant. vel. enhancement factor & {\small Effect modulated by $n$}\\
$\epsilon$ & Ret. vel. enhancement factor & {\small Effect modulated by $m$} \\
$\lambda_1$ & Diffusivity barrier, AIS & {\small Reduces the rate of diffusion within the axon initial segment} \\
$\lambda_2$ & Diffusivity barrier, SC & {\small Reduces the rate of diffusion within the synaptic cleft} \\
 \hline
\end{tabular}
\caption{Glossary of symbols for the single-edge model, adapted from~\cite{Torok2021}. 
Values marked with an asterisk were estimated by~\cite{Konzack2007} in cultured rodent neurons. 
Ant. = anterograde, ret. = retrograde, conc. = concentration, vel. = velocity.}
\label{Table_for_Variables}
\end{table}
\clearpage
\clearpage
\clearpage
\clearpage
\clearpage
\clearpage
\clearpage%
\clearpage
\begin{figure}
\centering
\includegraphics[width=1\linewidth]{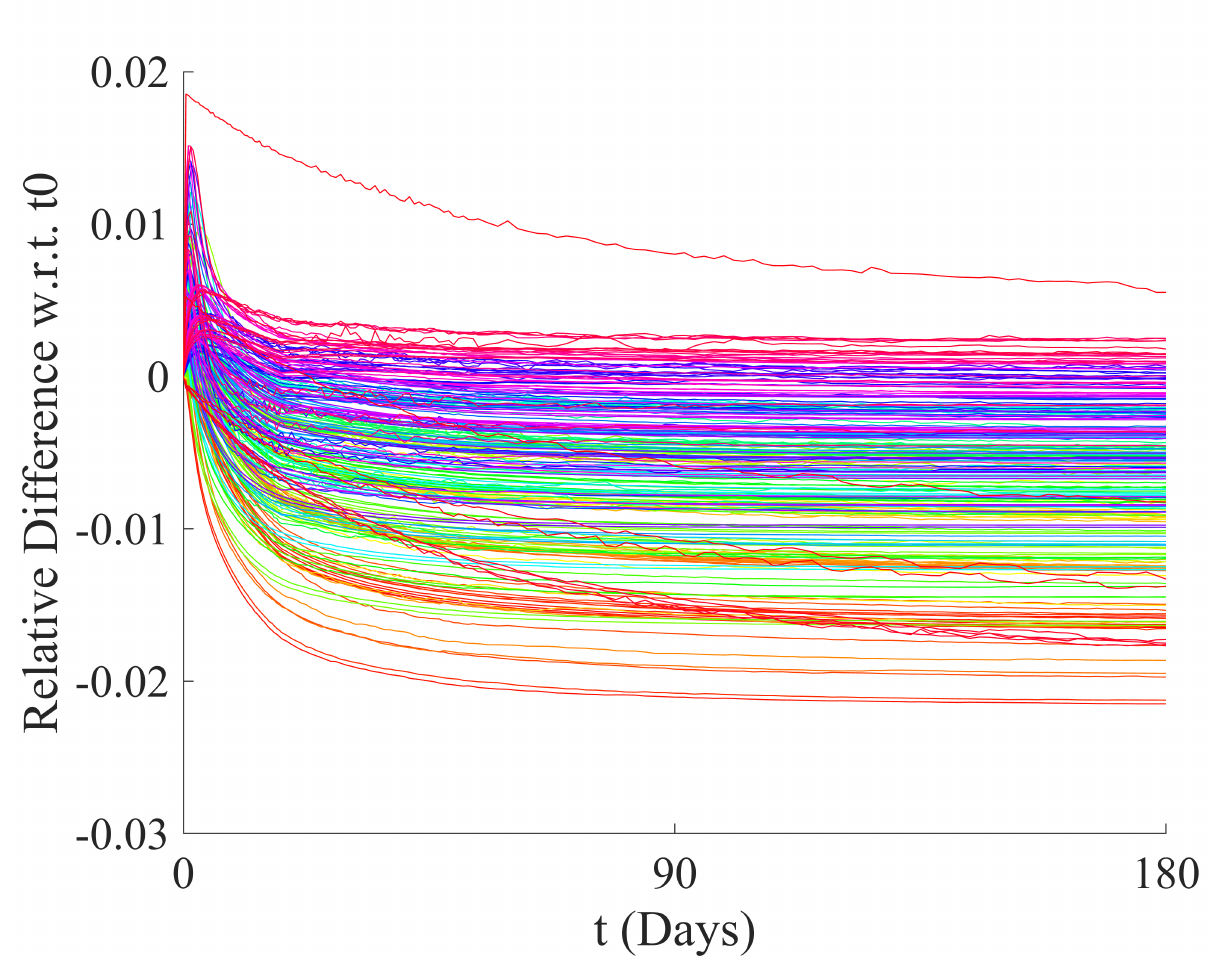}
\caption{Mass conservation of the NTM. Line plots of the relative error in the total mass of the system over time with respect to the initial total mass.}
\label{fig:masscons}
\end{figure}
\clearpage

\appendix
\section{The discretized problem}
Following \cite{Torok2021} we introduce an in-homogeneous  spatial mesh on each edge comportment of the form: 
$x_{\ell}=x_{K_{\ell}}<x_{K_{\ell}+1}<\dots<x_{K_{\ell+1}-1}<x_{K_{\ell+1}}=x_{\ell+1},\text{for}\;\; \ell=0, \dots,3$ and $x_{K_4}<\dots<x_{K-1}<x_K=L$ . In particular:
$x_0<\dots<x_{K_1}=x_1<\dots<x_{K-1}<x_K=L$\;.\\
We then introduce a time grid $t_p=p \Delta t, p=0,1,2,\dots$\\
Approximations $N_i^p, M_i^p, i=1, \dots h $ to the solutions $N_i(t_p), M_i(t_p), i=1 \dots h$ of \eqref{quasi-static model} are obtained by means of a 
first order Euler method leading to:
\begin{equation}
\begin{cases}
N_{i}^{p+1}=N_i^p+\frac{\sum_j \left(c_{ji}J_{ji}^p-c_{ij}J_{ij}^p\right)}{\left(V_i \left(1+\frac {\gamma_1 N_{i}^p(2\beta-\gamma_2 N_{i}^p)}{(\beta-\gamma_2 N_{i}^p)^2}\right)+
\sum_j \left(C_{ij}^{i,p}+C_{ji}^{i,p}\right)\right)} \Delta t\\
M_i^p=\frac {\gamma_1 (N_i^p)^2}{\beta-\gamma_2 N_i^p} \\
N_i(0)=N_{0i}\,.
\end{cases}
\end{equation}
%where the superscript $^p$ indicates the solution at time $t_p$ and 
where $C_{ij}^{i,p} \approx C_{ij}^i(t_p), C_{ji}^{i,p} \approx C_{ji}^i(t_p), J_{ij}^p \approx J_{ij}(t_p)$ are approximations of \eqref{Cij}, \eqref{Cji} and \eqref{constant flux} at time $t_p$.  In addition, 
$C_{ij}^{i,p}, C_{ji}^{i,p}, J_{ij}^p$ are determined by the approximated solutions $(n_{ij}^p, m_{ij}^p)$  of \eqref{equ n, vanishing rho} with Dirichlet boundary conditions $n_{ij,0}^p=N_i^p$ and $n_{ij,L}^p=N_j^p$. 
Furthermore, denoting by $q^{i,p}_{ij}, q^{i,p}_{ji}$  the approximated solutions of \eqref{q_new1} determined respectively by the approximated steady states solutions $n_{ij}=n_{ij}^p,n_{ji}=n_{ji}^p$ and 
boundary condition equal to $1$ at vertex $P_i$ and $0$ at vertex $P_j$,  
a trapezoidal integration method leads to:
\begin{align*}
C_{ji}^{i,p}=
%\frac{1}{2} \sum_{k=1}^{K-1}(x_k-x_{k-1} )(q^{i,n}_{ji}(x_k) +q^{i,n}_{ji}(x_{k+1})+\frac{1}{2}  
&\frac{c_{ji}}{2}\sum_{\scriptstyle (0\leq k \leq K_3-1)\cup(K_4\leq k \leq K_L-1)}(x_k-x_{k-1}) \big[q^{i,p}_{ji,k-1}\left(1+\frac{\gamma_1 n_{ji,k-1}^p(2\beta-\gamma_2  n_{ji,k-1}^p)}{(\beta-\gamma_2  n_{ji,k-1}^p)^2}\right)
+\\
&c_{ji} q^{i,p}_{ji,k}\left(1+\frac{\gamma_1 n_{ji,k}^p(2\beta-\gamma_2  n_{ji,k}^p)}{(\beta-\gamma_2  n_{ji,k}^p)^2}\right)
 \big]+ \frac{c_{ji}}{2} \sum_{K_3\leq k \leq K_4-1}(x_k-x_{k-1} )(q^{i,p}_{ji,k-1} +q^{i,p}_{ji,k})
 \end{align*}
and
\begin{align*}
C_{ij}^{i,p}=&\frac{c_{ij}}{2}\sum_{\scriptstyle (0\leq k \leq K_3-1)\cup(K_4\leq k \leq K_L-1)}(x_k-x_{k-1})\big[
q^{i,p}_{ij,k-1} \left(1+\frac{\gamma_1 n_{ij,k-1}^{p}(2\beta-\gamma_2  n_{ij,k-1}^p)}{(\beta-\gamma_2  n_{ij,k-1}^p)^2}
 \right)+\\
 & c_{ij} q^{i,p}_{ij,k} \left(1+\frac{\gamma_1 n_{ij,k}^{p}(2\beta-\gamma_2  n_{ij,k}^p)}{(\beta-\gamma_2  n_{ij,k}^p)^2}
 \right) \big]+ \frac{c_{ij}}{2} \sum_{K_3\leq k \leq K_4-1}(x_k-x_{k-1} )(q^{i,p}_{ij,k-1} +q^{i,p}_{ij,k})
 \end{align*}
where $n_{ij,k}^p \approx n_{ij}(x_k,t_p)$ and  $q_{ij,k}^{i,p} \approx q_{ij}^i(x_k,t_p)$.
%In particular, $n_{ij}^p$ is the solution of \eqref{steady_states_n} with Dirichlet boundary conditions $n_{ij,0}^p=N_i^p$ and $n_{ij,L}^p=N_j^p$, whereas $A_{ij}^p$ is its constant flux.
Setting $\Delta x_{k+1}=x_{k+1}-x_k$ we see that:
\begin{align}\label{n_approx}
n_{ij,k+1}^p  \cong n_{ij,k}^p+\frac{\Delta x_{k+1}}{a(x_k)}\left( -h(x_k,n_{ij,k}^p)-J_{ij}^p\right)
\end{align}
and
\begin{align}\label{q_approx}
q_{ij,k+1}^{i,p} \cong q_{ij,k}^{i,p}+\frac{\Delta x_{k+1}}{a(x_k)}\left( -h_n(x_k,n_{ij,k}^p)q_{ij,k}^{i,p}-W_{ij}^{i,p}\right)
\end{align}
where $W_{ij}^{i,p}$ is the approximation of the integration constant of \eqref{q_new1} determined by the approximated steady state $n_{ij}^p$
and boundary conditions equal to $1$ at vertex $P_i$ and $0$ at vertex $P_j$. 

%{\color{blue} In the numerical %implementation, in order to compute the approximated solutions $n_{ij}^p, q_{ij}^p$ we make use of the matlab solver ode 45 which is based on an explicit Runge-Kutta (4,5) formula giving a more accurate approximation of $n_{ij,k}^p, n_{ij,k}^p$ than \eqref{n_approx}, \eqref{q_approx}.} {\color{red} Do we need to specify this?}
%\\
%Finally, from \eqref{eq q_ji}, \eqref{eq d_A_ij} we see that
%\begin{align*}
% J_{ij}^{p+1} =J_{ij}^p+ (N_i^{p +1}-N_i^p)W_{ij}^{i,p} + (N_j^{p +1}-N_j^p)W_{ij}^{j,p}\;.
 %\end{align*}

%{\color{red} Do we insert this formula even if in the numerical simulation we do the shooting procedure at each time step to calculate $A_{ij}^p$ since it provides a more accurate approximation?}
%{\color{blue} If we insert it but don't use it, we must at least mention that we don't use it.}

%\section{Equations for the ``continuous in time'' model}
%{\bf Equations for $m,n\ge 0$ in the nodes.} In node $i$, 
%under the hypothesis of quasi steady states also in the nodes:
%$$
%(m_i)_t=-\Gamma(m_i,n_i)=-\beta m_i+ n_i(\gamma_1n_i + \gamma_2m_i)=m_i(\gamma_2 n_i -\beta)+\gamma_1 n_i^2
%$$
%\begin{equation}\label{n equ}
%(n_i)_t=\,\text{diffusion term}\,+\Gamma(m_i,n_i), %=-(\text{flux})_x-m(\gamma n -\beta)-\gamma n^2
%\end{equation}
$$
%\Gamma(M_i,N_i)=\beta M_i- N_i(\gamma_1N_i + \gamma_2M_i)=-M_i(\gamma_2 N_i -\beta)-\gamma_1 N_i^2=0
$$

\end{document}